\documentclass[12pt, notitlepage]{amsart}
\usepackage{latexsym, amsfonts, amsmath, amssymb, amsthm, cite}

\newtheorem{theorem}{Theorem}

\usepackage{graphicx}
\usepackage{mathrsfs}

\begin{document}

\title{Bertrand's postulate and the existence of finite fields} 
 \author{K. Soundararajan} 
\address{Department of Mathematics \\ Stanford University \\
450 Jane Stanford Way, Bldg. 380\\ Stanford, CA 94305-2125}
\email{ksound@stanford.edu}
\thanks{I am partially supported through a grant from the National Science Foundation and a Simons Investigator award from the Simons Foundation.} 
  
 \maketitle
 
 While the finite field ${\Bbb F}_p = {\Bbb Z}/p{\Bbb Z}$ is  easily constructed (here $p$ denotes  a prime number), it is less easy to see the  existence of finite fields of order $p^n$ for higher prime powers.  One way  to construct such fields is to consider the polynomial ring ${\Bbb F}_p[x]$, which is readily seen to possess a division algorithm, and is thus a Euclidean domain and therefore a unique factorization domain.  If we can find an irreducible polynomial $P \in {\Bbb F}_p[x]$ of degree $n$, then the quotient ring ${\Bbb F}_p[x]/(P)$ gives a finite field of size $p^n$.   Thus, in this approach, the existence of finite fields of prime power order amounts to the existence of an irreducible polynomial in ${\Bbb F}_p[x]$ of any given degree $n$.

The problem of finding an irreducible polynomial $a_n x^n+\ldots+a_0 \in {\Bbb F}_p[x]$ may be viewed as analogous to finding a prime number of the form $a_n p^n + \ldots +a_0$ with each $a_i \in \{0, 1, \ldots, p-1\}$; in other words, a prime number with $n$ digits in its base $p$ expansion.    This problem for prime numbers is hardest in base $2$, and asks for a prime number with a given number (at least $2$) of binary digits.    The existence of such primes is guaranteed by 
 Bertrand's postulate (first established by Chebyshev) which asserts that for every natural number $N \ge 2$, there is a prime number between $N$ and $2N$.   A lovely proof of Bertrand's postulate was found by Ramanujan, and further simplified by Erd{\H o}s, and a lucid treatment of this ``book proof" may be found in \cite{AZ}.    Motivated by the Erd{\H o}s--Ramanujan proof, this note gives a quick proof of the existence  of an irreducible polynomial in ${\Bbb F}_p[x]$  of degree $n$.   We first present the proof, and then discuss its  parallels with Bertrand's postulate.   In discussing irreducible polynomials, we may naturally restrict attention to monic polynomials (that is, polynomials with leading coefficient $1$).  

\begin{theorem}  \label{thm1} Let ${\Bbb F}$ be a finite field.  For any natural number $n$, there exists a monic irreducible polynomial of degree $n$ in ${\Bbb F}[x]$.   
\end{theorem} 

Starting with the field ${\Bbb F}_p = {\Bbb Z}/p{\Bbb Z}$, we may find a monic irreducible polynomial $P \in {\Bbb F}_p[x]$ of degree $n$.  As mentioned already, this shows the existence of a finite field of size $p^n$ by considering the quotient ${\Bbb F}_p[x]/(P)$.  

\begin{proof}   Let $q$ denote the size of the field ${\Bbb F}$.   For any natural number $n$, define ${\mathcal F}(n)$ to be the product of all monic polynomials of degree $n$ in 
${\Bbb F}[x]$.  Since there are $q^n$ monic polynomials of degree $n$, ${\mathcal F}(n)$ is a polynomial of degree $n q^n$.  We claim that 
\begin{equation} 
\label{1} 
\frac{{\mathcal F}(n)}{{\mathcal F}(n-1)^q} = \prod_{\substack{P \\ d=\text{deg}(P) |n}} P,
\end{equation} 
where the product is over monic irreducible polynomials $P$ such that $\text{deg}(P)$ divides $n$.  

Granting the claim, we now establish the theorem.   The case $n=1$ is evident: there are $q$ monic polynomials of degree $1$ and all are irreducible.   Suppose now that $n\ge 2$.   The claimed relation shows that ${\mathcal F}(n)/{\mathcal F}(n-1)^q$ is in fact a polynomial in ${\Bbb F}[x]$.   The degree of the left side of \eqref{1} is $nq^n - q((n-1)q^{n-1}) = q^n$.  Suppose that there are no monic irreducibles of degree $n$.  Then the right side of \eqref{1} includes only monic irreducibles with degree $\le \lfloor n/2\rfloor$, since a divisor $d$ of $n$ with $d<n$ must satisfy $d\le \lfloor n/2\rfloor$.   Thus the right side of \eqref{1} divides the product over all monic polynomials of degree $\le \lfloor n/2\rfloor$.  Since there are $q^d$ monic polynomials of degree $d$, we conclude that the degree of the right side of \eqref{1} must be 
$$ 
\le \sum_{1\le d\le \lfloor n/2\rfloor} dq^d \le \lfloor n/2 \rfloor \sum_{1\le d\le \lfloor n/2\rfloor} q^d = \lfloor n/2 \rfloor 
\frac{q^{\lfloor n/2\rfloor +1}-q}{q-1} < \lfloor n/2 \rfloor \frac{q}{q-1} q^{\lfloor n/2\rfloor} \le 2 \lfloor n/2\rfloor q^{\lfloor n/2\rfloor}.
$$ 
Thus, if there are no irreducibles of degree $n$, we must have 
$$ 
q^n < 2 \lfloor n/2 \rfloor q^{\lfloor n/2 \rfloor}. 
$$ 
However $q^{n} \ge q^{\lfloor n/2\rfloor} 2^{\lfloor n/2 \rfloor} \ge 2 \lfloor n/2 \rfloor q^{\lfloor n/2\rfloor}$ (since $2^k \ge 2k$ for all $k\ge 1$), which is a contradiction.   Therefore there must be a monic irreducible of degree $n$.


We now prove the claimed relation \eqref{1}.  To this end, let $P$ be a monic irreducible of degree $d$ (not necessarily a divisor of $n$), and consider the exact power of $P$ that divides ${\mathcal F}(n)$.  This exponent equals the number of monic polynomials of degree $n$ that are divisible by $P$, plus the number of such polynomials divisible by $P^2$, plus the number divisible by $P^3$, and so on.  Now, a monic polynomial $F$ of degree $n$ is divisible by $P^k$ if we can write $F=P^k G$ for a monic polynomial $G$ of degree $n-kd$ (assuming this is non-negative), and there are $q^{n-kd}$ such $G$.  Thus the number of monic polynomials of degree $n$ divisible by $P^k$ may be expressed as $\lfloor q^{n-kd}\rfloor$ --- when $n-kd <0$ one obtains the floor of a quantity below $1$, which is zero.  We conclude that the exact power of $P$ dividing ${\mathcal F}(n)$ equals 
\begin{equation} 
\label{2} 
\sum_{k=1}^{\infty} \lfloor q^{n-kd} \rfloor.  
\end{equation} 

Applying \eqref{2} to determine the power of $P$ dividing ${\mathcal F}(n)$ and ${\mathcal F}(n-1)^q$ and subtracting, we find that the power of $P$ dividing ${\mathcal F}(n)/{\mathcal F}(n-1)^q$ equals 
\begin{equation*} 
\sum_{k=1}^{\infty} \Big( \lfloor q^{n-kd} \rfloor - q \lfloor q^{n-1 - kd} \rfloor \Big). 
\end{equation*} 
If $kd \le n-1$, or if $kd \ge n+1$ then $\lfloor q^{n-kd} \rfloor - q \lfloor q^{n-1 - kd} \rfloor  =0$, while if $kd =n$ (which occurs when $d$ is a divisor of $n$) $\lfloor q^{n-kd} \rfloor - q \lfloor q^{n-1 - kd} \rfloor  =1$.  This establishes the claim \eqref{1}, and hence the theorem. 
\end{proof}

 We now point out the parallels between the proof given above and the Erd{\H os}--Ramanujan proof of Bertrand's postulate.  
 The key idea there is to consider the middle binomial coefficient $\binom{2N}{N}$, and to compare its size with its prime factorization.  Since the middle binomial coefficient is the largest entry in its row of Pascal's triangle, one sees that 
 $2^{2N}/(2N+1) \le \binom{2N}{N} \le 2^{2N}$, which gives good control on its size.  To understand the prime factorization of $\binom{2N}{N}$, we determine the exact power of a prime $p$ that divides $n!$.   Since there are $\lfloor n/p \rfloor$ multiples of $p$ below $n$, $\lfloor n/p^2 \rfloor$ multiples of $p^2$, and so on, it follows that the exact power of $p$ dividing $n!$ is 
 \begin{equation} 
 \label{4}
\Big \lfloor \frac np\Big \rfloor +\Big \lfloor\frac{n}{p^2} \Big\rfloor + \ldots = \sum_{k=1}^{\infty}\Big \lfloor \frac{n}{p^k} \Big\rfloor.
 \end{equation} 
 Applying \eqref{4} with $n=2N$ and $n=N$, we may determine the prime factorization of $\binom{2N}{N}$.   If there are no primes between $N$ and $2N$, then from this prime factorization we may extract an upper bound for $\binom{2N}{N}$ which contradicts the lower bound $\binom{2N}{N} \ge 4^N/(2N+1)$.  This is a bare bones description of the proof of Bertrand's postulate, and we refer to \cite{AZ} for a fleshed out account.

Revisiting the argument of Theorem \ref{thm1}, the polynomial ${\mathcal F}(n)$ is meant to resemble a factorial --- the analogy is not exact but only approximate, since we multiply only the monic polynomials of degree exactly $n$, rather than degree at most $n$.   The quotient ${\mathcal F}(n)/{\mathcal F}(n-1)^q$ plays the role of the binomial coefficient $\binom{2N}{N}$, and the analogue of the size of $\binom{2N}{N}$ is the degree of ${\mathcal F}(n)/{\mathcal F}(n-1)^q$.  While factorials grow faster than exponentially, the binomial coefficient $\binom{2N}{N}$ grows only exponentially.   Similarly, while ${\mathcal F}(n)$ has degree $nq^{n}$, the degree of ${\mathcal F}(n)/{\mathcal F}(n-1)^q$ is the slower growing quantity $q^n$.  Finally, our proof plays off the degree of ${\mathcal F}(n)/{\mathcal F}(n-1)^q$ with its factorization into irreducibles \eqref{1}.  This is achieved by finding the exact power of an irreducible $P$ dividing ${\mathcal F}(n)$, which is computed in \eqref{2} mirroring the computation in the case of integers in \eqref{4}.

This note arose through my experiences with teaching Math 62DM at Stanford, which is the middle course of a three Quarters long honors sequence aimed at incoming undergraduates with a strong interest in mathematics.   The DM sequence focuses on ``Discrete Methods" in modern mathematics, and this particular course is centered around finite fields and their applications in combinatorics and elsewhere.  I begin with a discussion of rings and factorization into irreducibles, concentrating on the Euclidean domains ${\Bbb Z}$ and polynomials over a field.  This leads naturally to preliminary investigations on primes in the integers, culminating in the Erd{\H o}s--Ramanujan proof of Bertrand's postulate.  Congruences and more general quotient rings are next discussed, leading to the construction of finite fields as quotients of ${\Bbb F}_p[x]$ by irreducible polynomials.  At this juncture, one needs the existence of an irreducible polynomial of degree $n$ in ${\Bbb F}_p[x]$.   One standard proof is to consider the zeta function 
$$ 
\zeta(s,{\Bbb F}) = \sum_{F \text{ monic}} \frac{1}{q^{s \text{deg}(F)}} = \prod_{P  \text{ monic irreducible}} \Big( 1- \frac{1}{q^{s\text{deg}(P)}} \Big)^{-1}.
$$ 
Since there are $q^n$ monic polynomials of degree $n$, from the series representation we obtain $\zeta(s,{\Bbb F}) = (1-q^{1-s})^{-1}$.   A consideration of $\log \zeta(s,{\Bbb F})$ as a series over the irreducibles $P$ leads to Gauss's relation 
\begin{equation} 
\label{5} 
q^n = \sum_{d|n} d\pi(d, {\Bbb F}) 
\end{equation}  
where $\pi(d,{\Bbb F})$ denotes the number of monic irreducibles of degree $d$.  Without a proper discussion of convergence of series and products, students found this argument difficult, and a  formal approach is not entirely accurate since one must take logarithms.   While setting an exam for the course, it occurred to me that the earlier proof of Bertrand's postulate may be adapted as explained above, and that the parallels might illuminate both results.  

As further pedagogical points, let us note that the relation \eqref{5} follows at once by comparing the degrees of both sides of \eqref{1}.   A discussion of M{\" o}bius inversion and determining a formula for $\pi(n, {\Bbb F})$ is naturally suggested.  Finally, after a more detailed discussion of finite fields, students will recognize ${\mathcal F}(n)/{\mathcal F}(n-1)^q$ as the polynomial $x^{q^n} -x$.



\begin{thebibliography}{20} 
 
 
 \bibitem{AZ}  M. Aigner and G.M.  Ziegler.  {\sl Proofs from the Book. }  Fourth edition.  Springer-Verlag, Berlin, (2010).
  
  
 
 \end{thebibliography}

 \end{document}